\documentclass[a4paper]{article}
\usepackage{theorem}
\usepackage[fleqn]{amsmath}
\usepackage{amsfonts,amssymb,amsxtra} 
\usepackage{hyperref}

\hypersetup{
	pdftitle = {The Kodaira Dimension of Certain Moduli Spaces of Abelian Surfaces},
	pdfauthor = {C. Erdenberger <erdenber@math.uni-hannover.de>},
	pdfsubject = {MSC(2000) 14K10},
	pdfkeywords = {abelian surfaces, moduli space, Kodaira dimension, general type}	 
}

\theoremstyle{plain}
\newtheorem{Def}{Definition}
\newtheorem{Cor}[Def]{Corollary}
\newtheorem{Lem}[Def]{Lemma}
\newtheorem{Prop}[Def]{Proposition}
\newtheorem{Thm}[Def]{Theorem}

\newenvironment{proof}{\textit{Proof. }}{\hspace*{\fill}{$\square$}\bigskip}

\renewcommand{\Im}{\mathrm{Im}}
\renewcommand{\mod}[1]{\; ( \bmod \; #1 )}
\newcommand{\wt}{\widetilde}

\newcommand{\bbC}{\mathbb{C}}
\newcommand{\bbH}{\mathbb{H}}
\newcommand{\bbN}{\mathbb{N}}
\newcommand{\bbQ}{\mathbb{Q}}
\newcommand{\bbR}{\mathbb{R}}
\newcommand{\bbZ}{\mathbb{Z}}

\newcommand{\sieg}[1][2]{\bbH_#1}

\newcommand{\Apol}[1]{\mathcal{A}_{#1}}
\newcommand{\Ap}{\Apol{p}}
\newcommand{\Ga}[1]{\Gamma_{#1}}
\newcommand{\Gp}{\Ga{p}}
\newcommand{\idm}{\mathbf{1}}
\newcommand{\McuspGp}[1][n]{\mathcal{M}^\ast_{#1}(\Gp)}
\newcommand{\otn}{\omega^{\otimes n}}
\newcommand{\SigL}{\Sigma_L}
\newcommand{\tSigL}{\wt{\Sigma}_L}
\newcommand{\Apk}{\Ap^\ast}
\newcommand{\tApk}{\wt{\mathcal{A}}_{\mathnormal{p}}^\ast}
\newcommand{\bcD}[1]{D^\circ(#1)}
\newcommand{\loc}{\mathcal{O}}
\newcommand{\hum}{\mathcal{H}}
\newcommand{\thumk}{\wt{H}}
\newcommand{\humC}[1]{\mathcal{C}_{#1,1}}
\newcommand{\humCk}[1]{\wt{C}_{#1,1}}
\newcommand{\resptot}{\bar{\mathcal{A}}_p^\ast}
\newcommand{\TSigL}{T_{\tSigL}}
\newcommand{\TtSigL}{T_{\tSigL}}
\newcommand{\GGp}{\mathcal{G}_p}

\newcommand{\orb}{\mathrm{orb}}
\newcommand{\id}{\mathrm{id}}

\newcommand{\Sp}{\mathrm{Sp}}
\newcommand{\Stab}{\mathrm{Stab}}
\newcommand{\Sym}{\mathrm{Sym}}
\newcommand{\Fix}{\mathrm{Fix}}

\begin{document}

\title{The Kodaira Dimension of Certain Moduli Spaces of Abelian Surfaces}

\author{C. Erdenberger}

\maketitle

\section{Introduction}

Moduli Spaces for polarized abelian surfaces are obtained as quotients of the Siegel upper half--plane $\sieg$ by 
arithmetic subgroups of the symplectic group $\Sp(4,\bbQ)$. In this paper we will look at the Kodaira dimension of the
moduli spaces $\Apol{t}$ which parameterize abelian surfaces with a polarization of type $(1,t)$. For many small values
of $t$ it is known that $\Apol{t}$ is rational or unirational (e.g. see the result by Gross and Popescu \cite{gp:epas98}).
On the other hand, there are lower bounds for $t$ such that $\Apol{t}$ is of general type. For instance, O'Grady 
\cite{og:kdmsas89} showed in 1989 that in the case $t = p^2$, $p$ a prime, the moduli space $\Apol{p^2}$ is of general 
type for $p \ge 17$ (which was later improved to $p \ge 11$ by Gritsenko and Sankaran \cite{gs:masp97}).

We will concentrate on the case where $t = p$ is a prime. There is a result due to Sankaran \cite{san:mpas97} which proves
that $\Ap$ is of general type for $p \ge 173$. While it is known that for small $p$ the spaces are rational or unirational,
there is a certain range of values where little is known. We will fill most of this gap by improving Sankaran's bound to 
$p > 71$. Moreover, we will show for some smaller values of $p$ that $\Ap$ is also of general type.

The structure of our proof is essentially the same as in Sankaran's article. We will, however, use a different
compactification of $\Ap$ which will enable us to get rid of the obstructions coming from the singularities in the
boundary. Our main tools will be Sankaran's results and the description of the singularities in $\Ap$ studied in 
\cite{bra:smafpt94} by Brasch.

This paper is organized as follows. In Section~\ref{sec:prelim} we establish notation and give an outline of the proof.
A toroidal compactification of $\Ap$ is constructed in Section~\ref{sec:torcomp} as a first step to finding a smooth 
projective model. In Section~\ref{sec:extbd} we apply a method due to Gritsenko that enables us to deal with the smooth
points in the boundary of the compactification. Section~\ref{sec:singint} calculates the obstructions arising inside $\Ap$.
The singularities in the boundary are dealt with in Section~\ref{sec:singbd}. Finally, in Section~\ref{sec:concl} we 
assemble the parts of the proof and establish our main result.

At this point I would like to thank Prof. Klaus Hulek for his continuous interest and support and all the others who
contributed to this paper.

\section{Preliminaries}	\label{sec:prelim}

Let $\sieg$ denote the Siegel upper half--plane of degree 2,
\begin{equation*}
	\sieg:= \left\{ \tau = \begin{pmatrix} \tau_1 & \tau_2 \\ \tau_2 & \tau_3 \end{pmatrix}
	\in \Sym(2,\bbC) \, ; \; \Im(\tau) > 0 \right\} \; ,
\end{equation*}
and define for any prime $p$ the group $\Gp$ as the subgroup of the symplectic group $\Sp(4,\bbQ)$ given by
\begin{equation*}
	\Gp := \left\{ M \in \Sp(4,\bbQ) \, ; \; M \in
	\begin{pmatrix}
		\bbZ & \bbZ & \bbZ & p \bbZ \\
		p \bbZ & \bbZ & p \bbZ & p \bbZ \\
		\bbZ & \bbZ & \bbZ & p \bbZ \\
		\bbZ & \frac{1}{p}\bbZ & \bbZ & \bbZ
	\end{pmatrix}
	\right\} \; .
\end{equation*}
The action of $\Gp$ on $\sieg$ given by
\begin{equation*}
	M = \begin{pmatrix} A & B \\ C & D \end{pmatrix} \, : \, \tau \mapsto (A \tau + B) (C \tau + D)^{-1}
\end{equation*}
is properly discontinuous and we denote the quotient $\sieg / \Gp$ by $\Ap$. The normal complex analytic space $\Ap$ is a
moduli space for $(1,p)$--polarized abelian surfaces. It has been studied in \cite{hw:bas89} and \cite{bra:smafpt94}.

We are concerned with determining those values of $p$ for which $\Ap$ is of general type. Towards this goal we will 
consider for each $n \in \bbN$ certain modular forms on $\sieg$ and try to extend these to $n$--canonical forms on a
smooth projective model of $\Ap$. By studying the behavior of the space of forms which can be extended in this way with
respect to $n$ we obtain a lower bound for the Kodaira dimension of $\Ap$. If this bound turns out to be maximal we can
conclude that $\Ap$ is of general type.

Let $\omega := d\tau_1 \wedge d\tau_2 \wedge d\tau_3$ denote the usual differential $3$--form on $\sieg$. If $f$ is a 
weight~$3n$ cusp form for $\Gp$ then $f\otn$ is a differential form on $\sieg$ which is invariant under the operation of
$\Gp$. Hence it gives an $n$--canonical form on $\Ap$ where the quotient map $\sieg \rightarrow \Ap$ is unbranched. In 
order to extend this to a resolution of singularities of a suitable compactification of $\Ap$, i.e. a smooth projective
model of $\Ap$, we have to overcome three obstacles. First, we need to be able to extend the form to the smooth points of
the boundary, second, over the branch locus and, third, over the singularities in the boundary. We will deal with these 
three problems in sections~\ref{sec:extbd} - \ref{sec:singbd} respectively. The first can be solved by choosing special
cusp forms in the manner of \cite{san:mpas97}. For the second we can proceed similar to \cite[Section~5]{san:mpas97}.
The last one will be dealt with by choosing a "good" compactification of $\Ap$.

We will now introduce the cusp forms we want to work with. For this we need a nontrivial cusp form of weight $2$ for 
$\Gp$. The following proposition deals with the existence of such a form.

\begin{Prop}	\label{prop:nontrivcusp}
	There is a nontrivial weight~$2$ cusp form for $\Gp$ if $p > 71$ or $p=37, 43, 53, 61$ or $67$.
\end{Prop}
\begin{proof}
	The proof for $p > 71$ is given in \cite[Proposition~3.1]{san:mpas97}. It is based on a result by Gritsenko 
	\cite[Theorem~3]{gri:imspas94} which shows that a Jacobi cusp form of weight $2$ and index $p$ can be lifted
	to a weight~$2$ cusp form for $\Gp$. Using an explicit formula for the dimension of the space of these Jacobi cusp
	forms and checking that it is positive for $p > 71$, Sankaran proves the existence of the desired nontrivial 
	weight~$2$ cusp form for $\Gp$. However, it is easy to check that this is also the case for the given smaller 
	values of $p$.
\end{proof}

Let $f_2$ be a nontrivial weight~$2$ cusp form for $\Gp$ and $f_n$ a cusp form of weight $n$. Then $f := f_2^n f_n$ is a 
cusp form of weight $3n$. The dimension of the space of such forms computes as the dimension of $\McuspGp$, the space of 
weight~$n$ cusp forms for $\Gp$, which was calculated by Sankaran.

\begin{Prop}	\label{prop:dimcusp}
	For $p \ge 3$ the space $\McuspGp$ of weight~$n$ cusp forms for $\Gp$ satisfies
	\begin{equation*}
		\dim \McuspGp = \frac{p^2 + 1}{8640} n^3 + O(n^2) \; .
	\end{equation*}
\end{Prop}
\begin{proof}
	\cite[Proposition~3.2]{san:mpas97}
\end{proof}

For the rest of this paper, $f$ will denote a weight~$3n$ cusp form as constructed above and the prime $p$ will be 
assumed to be greater than $3$.

\section{Toroidal compactification \texorpdfstring{of $\Ap$}{}} \label{sec:torcomp}

To determine the Kodaira dimension of $\Ap$ we have to find a smooth projective model of $\Ap$. In this section we will
as a first step construct a projective model by using the method of toroidal compactification which is due to Hirzebruch,
Mumford et al. (cf. \cite{amrt:sc75}). For a more detailed description of this method in the case of moduli spaces we 
refer the reader to \cite{hkw:ms93}.

The Tits building of $\Gp$ is represented by the following graph:
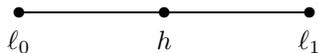
\begin{center}
	\begin{picture}(120,20)(0,0)
		\put(5,15){\circle*{4}}
		\put(5,8){\makebox(0,0)[t]{$\ell_0$}}
		\put(60,15){\circle*{4}}
		\put(60,8){\makebox(0,0)[t]{$h$}}
		\put(115,15){\circle*{4}}
		\put(115,8){\makebox(0,0)[t]{$\ell_1$}}
		\put(5,15){\line(1,0){110}}
	\end{picture}
\end{center}
The three vertices represent two $1$--dimensional isotropic subspaces $\ell_0$ and $\ell_1$ and one $2$--dimensional
isotropic subspace $h$ corresponding to proper rational boundary components of $\sieg$. Note that it is sufficient to 
construct an admissible fan for the corank~$2$ boundary component corresponding to $h$, since such a fan uniquely 
determines an admissible collection of fans needed for the compactification process. The Legendre fan is the standard 
choice for this (and the one that Sankaran chooses). It gives a toroidal compactification of $\Ap$ denoted by $\Apk$ in 
Sankaran's paper \cite{san:mpas97}.

An analysis of the singularities of $\Apk$ can be found in the thesis of Brasch \cite{bra:smafpt94} (see also 
\cite{san:mpas97}). It turns out that all the non--canonical singularities in the open corank~$2$ boundary component 
of $\Apk$ are coming from the fact that the Legendre fan is singular with respect to the lattice induced by the group 
$\Gp$. By choosing a suitable nonsingular subdivision of the Legendre fan one can therefore hope to be left only with 
mild (i.e. canonical) singularities in that component, so that the extension of pluricanonical forms will pose no further
problems there. As we will see in Section~\ref{sec:singbd} this is indeed the case. The following lemma gives us the 
existence of such a nonsingular subdivision.

\begin{Lem}
	There is a nonsingular subdivision $\tSigL$ of the Legendre fan $\SigL$ such that the toroidal compactification
	of $\Ap$ determined by $\tSigL$ is projective.
\end{Lem}
\begin{proof}
	We can apply \cite[Theorem~7.20]{nam:tcss80} to the Legendre fan $\SigL$ to obtain a nonsingular subdivision 
	$\smash{\tSigL}$ of $\SigL$ which is admissible with respect to $\Gp$. The proof of this theorem shows that the toric 
	variety defined by $\smash{\tSigL}$ is obtained by blowing up the toric variety defined by $\SigL$ with respect to 
	some coherent sheaf of ideals (cf. \cite[p.~32]{kkmsd:te73}). Hence, since the toroidal compactification of $\Ap$ 
	determined by the Legendre fan $\Apk$ is projective, so is the one determined by $\smash{\tSigL}$.
\end{proof}

We will denote the toroidal compactification of $\Ap$ determined by such a nonsingular subdivision $\smash{\tSigL}$ by 
$\smash{\tApk}$. There is a stratification of $\smash{\tApk}$ into disjoint sets:
\begin{equation*}
	\tApk = \Ap \amalg \bcD{\ell_0} \amalg \bcD{\ell_1} \amalg \bcD{h} \; ,
\end{equation*}
where the strata $\bcD{\ell_0}$ and $\bcD{\ell_1}$ denote the open corank~$1$ boundary components determined by $\ell_0$
and $\ell_1$ respectively and $\bcD{h}$ denotes the open corank~$2$ boundary component determined by $h$ (cf. 
\cite[Definition~3.78 and Remark~3.79]{hkw:ms93}). Note that since there is only one admissible fan for each of the 
boundary components other than the one determined by $h$, the toroidal compactifications in the direction of these 
boundary components coincide. Hence $\Apk$ and $\smash{\tApk}$ only differ in the open corank~$2$ boundary component 
$\bcD{h}$. This will enable us to use the results of Sankaran when we are dealing with obstructions to extending our
differential forms coming from the other boundary components.

\section{Extension to the boundary}	\label{sec:extbd}

We will now extend the cusp form $f = f_2^n f_n$ we introduced in Section~\ref{sec:prelim} to the toroidal 
compactification $\tApk$ constructed in the previous section. The fact that we have chosen $f$ in such a way that it has
a high order of vanishing at infinity guarantees us the desired extendibility to the generic points of the boundary.

\begin{Lem}	\label{lem:ext}
	Let $f:=f_2^n f_n$ be a weight~$3n$ cusp form for $\Gp$, where $f_2$ is a nontrivial weight~$2$ cusp form and $f_n$
	a cusp form of weight $n$. Then the differential form $f\otn$ extends over the smooth part of $\tApk$ away from any
	ramification.
\end{Lem}
\begin{proof}
	In \cite[Proposition~4.1]{san:mpas97} Sankaran only needed to consider the two corank~$1$ boundary components of
	$\Apk$ in order to extend the form $f$ to the whole smooth part of the boundary of $\Apk$, because in that 
	compactification of $\Ap$ the corank~$2$ boundary component has codimension $2$, so the form $f$ can always be 
	extended over the smooth points there. The proof is essentially an application of 
	\cite[Chapter~4, Theorem~1]{amrt:sc75}.

	We can use Sankaran's result to conclude that the form $f$ extends over the desired parts of the open corank~$1$
	boundary components $\bcD{\ell_0}$ and $\bcD{\ell_1}$ of $\smash{\tApk}$, since these coincide in both 
	compactifications as we observed in the previous section. However, the corank~$2$ boundary component $\bcD{h}$ of 
	$\smash{\tApk}$ has codimension $1$, so there is something to check here. To apply \cite[Theorem~1]{amrt:sc75} to 
	$\bcD{h}$ the way Sankaran did it for the corank~$1$ boundary components, we have to check that certain coefficients
	in the Fourier expansion of the form $f$ near the boundary component determined by $h$ vanish, i.e the form $f$ has 
	to vanish of order $n$ there. We can proceed in the same way Sankaran did for the corank~$1$ boundary components and
	write the expansion of $f$ as a product of expansions of $f_2$ and  $f_n$. The coefficients of the expansion of $f$
	can then be expressed in terms of those of the other two expansions and the required condition is now easy to verify
	(cf. \cite[Proposition~3.1]{gs:masp97}).
\end{proof}

\section{Singularities in \texorpdfstring{$\Ap$}{the interior}}	\label{sec:singint}

We next want to extend our forms over the points in the interior of $\smash{\tApk}$, i.e. $\Ap$ itself. In the previous 
section we extended the forms over the general points of $\smash{\tApk}$, but since the group $\Gp$ does not act freely on
$\sieg$, we had to exclude the branch locus. This causes two problems. First, the forms we are using collect poles along 
the branch divisors. Hence they will not give us global sections in the bundle $\loc_{\tApk} (n K_{\tApk})$ as desired but 
in a bundle determined by a suitable linear combination of the canonical divisor with these branch divisors. The second 
problem are the singularities in $\smash{\tApk}$. We need to find a resolution of these singularities and have to extend 
our forms over this resolution. We will start by calculating the branch divisors first and deal with the singularities 
later.

Since $\smash{\tApk}$ is a normal projective variety it is smooth in codimension $1$. This means that the only branch 
divisors are coming from quasi--reflections in $\Gp$. According to \cite[Kapitel~2, Folgerung~2.2]{bra:smafpt94} it 
suffices to consider involutions, i.e. elements of order $2$ different from $-\idm$. Two involutions in $\Gp$ are given by
\begin{equation*}
	I_1 :=
		\begin{pmatrix}
			-1 & 0 & 0 & 0 \\
			0 & 1 & 0 & 0 \\
			0 & 0 & -1 & 0 \\
			0 & 0 & 0 & 1
		\end{pmatrix}
	\quad \text{ and } \quad
	I_2 :=
		\begin{pmatrix}
			-1 & -1 & 0 & 0 \\
			0 & 1 & 0 & 0 \\
			0 & 0 & -1 & 0 \\
			0 & 0 & -1 & 1 
		\end{pmatrix} \; .
\end{equation*}
Their fixed loci in $\sieg$ are computed in \cite[p.~236]{hkw:smscas91} and are given by
\begin{align*}
	& \hum_1 := \Fix(I_1) =
		\left\{ \begin{pmatrix} \tau_1 & 0 \\ 0 & \tau_3 \end{pmatrix} \in \sieg 
		\, ; \; \tau_1 , \tau_3 \in \sieg[1] \right\} \; , \\
	& \hum_2 := \Fix(I_2) =
		\left\{ \begin{pmatrix} \tau_1 & - \frac{1}{2} \tau_3 \\ - \frac{1}{2} \tau_3 & \tau_3 \end{pmatrix} \in \sieg
		\, ; \; \tau_1 , \tau_3 \in \sieg[1] \right\} \; .
\end{align*}

We denote the closures of the images of $\hum_i$ for $i=1,2$ under the natural projections in $\smash{\tApk}$ by 
$\smash{\thumk_i}$. Since according to \cite[Proposition~5.2]{san:mpas97} every involution in $\Gp$ is conjugate to 
either $\pm I_1$ or $\pm I_2$ and $-\idm$ acts trivially on $\sieg$, $\smash{\thumk_1}$ and $\smash{\thumk_2}$ are the 
only branch divisors coming from the map $\sieg \rightarrow \Ap \hookrightarrow \smash{\tApk}$.

There might be other branch divisors in the boundary coming from the toroidal compactification of $\Ap$ (there is an
induced action of suitable subgroups of $\Gp$ on each of the toric varieties determined by the fans used in the 
compactification). However, it is easy to check that this is not the case (using arguments for groups acting on toric
varieties). Hence, $\smash{\thumk_1}$ and $\smash{\thumk_2}$ are in fact the only branch divisors in $\smash{\tApk}$. It 
also follows that they intersect the boundary $\smash{\tApk} \setminus \Ap$ of $\smash{\tApk}$ in codimension $\ge 2$. We 
can now regard our forms as sections in a certain bundle as follows:

\begin{Prop}	\label{prop:ext}
	If $n \in \bbN$ is sufficiently divisible, 
	$n \left( K_{\tApk} + \frac{1}{2} \thumk_1 + \frac{1}{2} \thumk_2 \right)$ is Cartier and $f\otn$ determines an
	element of
	\begin{equation*}
		H^0 \left( \tApk , n \left( K_{\tApk} + \frac{1}{2} \thumk_1 + \frac{1}{2} \thumk_2 \right) \right) \; .
	\end{equation*}
\end{Prop}
\begin{proof}
	Above $\smash{\thumk_1}$ and $\smash{\thumk_2}$ the map $\sieg \rightarrow \Ap$ is ramified with index $2$, so 
	$f \otn$ collects poles of order $\frac{n}{2}$ along $\smash{\thumk_1}$ and $\smash{\thumk_2}$. Since these are 
	the only branch divisors we can conclude using Proposition~\ref{lem:ext} that $f \otn$  defines a global section 
	in the bundle $\smash{n ( K_{\tApk} + \frac{1}{2} \thumk_1 + \frac{1}{2} \thumk_2 )}$.
\end{proof}

We will need to calculate numerical conditions which guarantee us that a section in this bundle is in fact an element
of $n K_{\tApk}$ as required for the purpose of determining the Kodaira dimension. However, we still need to resolve
the singularities of $\smash{\tApk}$. The blow--ups we will make to achieve this will produce exceptional divisors which 
will again impose certain numerical conditions on us. Therefore we will now deal with the singularities lying in the 
interior of $\smash{\tApk}$, i.e in $\Ap$ itself, first and then calculate the obstructions coming from the branch 
divisors and the exceptional divisors simultaneously. The singularities in the boundary of $\smash{\tApk}$ will be dealt
with in the next section.

The singularities in $\Ap$ are independent of the toroidal compactification chosen, so we can make use of Sankaran's 
results. Pluricanonical forms can always be extended at the canonical singularities, so it suffices to consider the 
non--canonical ones. According to \cite[Proposition~5.5]{san:mpas97} all the non--canonical singularities of $\Ap$ are
lying in $\smash{\thumk_1}$. The singularities in $\smash{\thumk_1} \cap \Ap$ are described in 
\cite[Kapitel~2, Hilfssatz~2.25]{bra:smafpt94}. They are exactly the points of four curves, $\humC{3}$, $\humC{4}$,
$\humC{5}$ and $\humC{6}$, forming a square. We denote the closures of their images in $\smash{\tApk}$ by 
$\smash{\humCk{i}}$, $i= 3, \dots , 6$, respectively. The non--canonical singularities of $\Ap$ lying in 
$\smash{\thumk_1}$ are precisely the points of two of these curves, $\smash{\humCk{3}}$ and $\smash{\humCk{5}}$, 
by \cite[Proposition~5.6]{san:mpas97}.

Exactly as in \cite[Section~5]{san:mpas97} we obtain two blow--ups $\smash{\wt{\beta}_1}$ and $\smash{\wt{\beta}_2}$ of 
$\smash{\tApk}$ which resolve the singularities at the general points of the four curves $\smash{\humCk{i}}$, 
$i=3, \dots , 6$. The singularities at the points where two of the curves meet are all canonical. By composing 
$\smash{\wt{\beta}_1}$ with $\smash{\wt{\beta}_2}$ we thus obtain a blow--up $\smash{\wt{\beta}} : \resptot \rightarrow
\smash{\tApk}$ such that $\resptot$ has at most canonical singularities over $\Ap$. Note that the curves 
$\smash{\humCk{i}}$, $i=3, \dots , 6$, do not intersect the open corank~$2$ boundary component $\bcD{h}$, so that 
$\smash{\wt{\beta}}$ defines an isomorphism there.

The obstructions coming from the exceptional divisors of this blow--up and from the branch divisors $\smash{\thumk_1}$ and 
$\smash{\thumk_2}$ are the same as in \cite[Section~5]{san:mpas97} and are given by

\begin{Lem}	\label{lem:obsint}
	The obstructions to extending a form $f \otn$ to a pluricanonical form on a resolution of singularities of 
	$\tApk$ lying in $\Ap$ are bounded by
	\begin{equation*}
		\left( \frac{7}{54} - \frac{1}{3p} \right) n^3 + O (n^2)
	\end{equation*}
\end{Lem}
\begin{proof}
	Using Proposition~\ref{prop:ext} we obtain a rational differential form on $\resptot$ which might fail to be 
	regular only at the exceptional divisors of $\smash{\wt{\beta}}$ or the branch divisors $\smash{\thumk_1}$ and 
	$\smash{\thumk_2}$. The resolutions of $\smash{\tApk}$ and $\Apk$ considered here and in Sankaran's paper coincide
	outside the open corank~$2$ boundary component $\bcD{h}$. Hence by \cite[Theorem~5.10]{san:mpas97} the form 
	$f \otn$ is under the given conditions regular away from $\bcD{h}$. This includes in particular the exceptional 
	divisors and the smooth points of $\smash{\thumk_1}$ and $\smash{\thumk_2}$ lying outside $\bcD{h}$. However, on 
	$\bcD{h}$ the blow--up $\smash{\wt{\beta}}$ defines an isomorphism. Since $\smash{\thumk_1}$ and $\smash{\thumk_2}$
	intersect the boundary in codimension $\ge 2$ as we observed earlier, we can therefore conclude that $f \otn$ can 
	be extended over the smooth points also in this component.
\end{proof}

\section{Singularities in the boundary}	\label{sec:singbd}

The blow--ups we have made in the previous section affect also the singularities in the boundary which will be examined
in this section. However, it will turn out that they have changed them in such a way that there are only canonical 
singularities in the boundary and that our forms can be extended without further obstructions.

We will abuse the notation and denote the strict transforms of the open boundary components $\bcD{\ell_0}$, 
$\bcD{\ell_1}$ and $\bcD{h}$ in $\resptot$ under the blow--up $\smash{\wt{\beta}}$ by $\bcD{\ell_0}$, $\bcD{\ell_1}$ and 
$\bcD{h}$ respectively. We will first investigate the open corank~$1$ boundary components $\bcD{\ell_0}$ and 
$\bcD{\ell_1}$. Since the spaces $\Apk$ and $\smash{\tApk}$ coincide here, we can again make use of Sankaran's results.

\begin{Prop}	\label{prop:cor1bdc}
	The open corank~$1$ boundary components $\bcD{\ell_0}$ and $\bcD{\ell_1}$ of $\resptot$ contain at most canonical
	singularities.
\end{Prop}
\begin{proof}
	We consider the singularities in the open corank~$1$ boundary components of $\Apk$ (which coincide with those in 
	$\smash{\tApk}$). They are described in \cite[Satz~4.6, 4.7]{bra:smafpt94}. In each of $\bcD{\ell_0}$ and 
	$\bcD{\ell_1}$ there are exactly four singular points. Four of them are isolated, the other four are lying on the
	intersections of the curves $\smash{\humCk{i}}$, $i=3, \dots , 6$, with the boundary components $\bcD{\ell_0}$ and
	$\bcD{\ell_1}$. Sankaran shows in the proof of \cite[Theorem~6.1]{san:mpas97} that the four isolated singularities 
	are canonical and that the other four are of the same type as the ones at the general points of the curves they are 
	lying on, so that they are resolved by $\smash{\wt{\beta}}$. Hence the open boundary components $\bcD{\ell_0}$ and 
	$\bcD{\ell_1}$ of $\resptot$ contain only canonical singularities as claimed.
\end{proof}

In order to calculate the singularities in the open corank~$2$ boundary component $\bcD{h}$ of $\resptot$, it suffices
to consider the ones in the component $\bcD{h}$ of $\smash{\tApk}$, since $\smash{\wt{\beta}} : \resptot \rightarrow 
\smash{\tApk}$ is an isomorphism when restricted to $\bcD{h}$. The component $\bcD{h}$ of $\smash{\tApk}$ as constructed
in the toroidal compactification is given as the quotient of an open subset of the toric variety $\TtSigL$ determined by 
the fan $\smash{\tSigL}$ by the action of a certain group which we denote by $\GGp$ following the notation of Brasch (cf. 
\cite[p.~118]{bra:smafpt94}). Since the fan $\smash{\tSigL}$ is nonsingular, so is its toric variety $\TtSigL$. Hence any
singularities in the quotient are coming from fixed points of the action of $\GGp$ on $\TtSigL$. We will calculate 
these quotient singularities and show that they are all canonical, thus proving that the component $\smash{\bcD{h}}$ of 
$\resptot$ has only canonical singularities.

We can decompose $\TtSigL$ in terms of its orbits as follows
\begin{equation*}
	\TtSigL = \coprod_{\wt{\sigma} \in \tSigL} \orb(\wt{\sigma}) \; .
\end{equation*}
Using this we can determine the singularities in the quotient $\TtSigL / \GGp$ by looking at the stabilizers in 
$\GGp$ of the cones $\wt{\sigma} \in \smash{\tSigL}$ and their actions on the corresponding orbits. Brasch has calculated
these stabilizers in \cite[Kapitel~3, Hilfssatz~5.23]{bra:smafpt94} for the Legendre fan $\SigL$. By 
\cite[Corollary~1.17]{oda:cbag88} we obtain an equivariant map $\id_\ast : \TtSigL \rightarrow \TSigL$ between the two
toric varieties associated to these two fans which will allow us to use Brasch's result.

\begin{Prop}	\label{prop:stabcont}
	Let $\wt{\sigma} \in \tSigL$ and $\sigma \in \SigL$ such that $\wt{\sigma} \subset \sigma$. If $\sigma$ is minimal
	in $\SigL$ with this property, then
	\begin{equation*}
		\Stab_{\wt{\sigma}} = \left\{ g \in \GGp \, ; \; g \cdot \wt{\sigma} = \wt{\sigma} \right\} 
		\subset \Stab_\sigma \; .
	\end{equation*}
\end{Prop}
\begin{proof}
	This is straightforward to check using the equivariance of the map $\id_\ast$.
\end{proof}

Note that we do not need to calculate all singularities in the quotient $\TtSigL / \GGp$, but only those that are lying
in the open subset which corresponds to the open corank~$2$ boundary component $\bcD{h}$. It turns out that these are 
exactly the ones coming from nontrivial stabilizers of those cones in $\smash{\tSigL}$ which are contained in the 
$2$-- or $3$--dimensional (but not in the $1$--dimensional) cones of the Legendre fan $\SigL$ (the others are identified
with points in the other boundary components by the gluing map used in the toroidal compactification). But the stabilizers
of these cones in the Legendre fan are by \cite[Kapitel~3, Hilfssatz~5.23]{bra:smafpt94} either trivial or cyclic of order 
$2$ or $3$. So in view of Proposition~\ref{prop:stabcont} the stabilizers of the relevant cones in $\smash{\tSigL}$ have 
the same property. Since stabilizers of order $2$ either yield smooth points or canonical singularities, we only need to 
worry about possible stabilizers of order $3$. Fortunately, there is only one cone in the Legendre fan $\SigL$, which has a 
stabilizer of order $3$. Following Brasch's notation, we denote it by $\sigma_{\lambda_3}$, where $\lambda_3$ satisfies
$\lambda_3^2 + \lambda_3 + 1 \equiv 0 \mod{p}$. Its stabilizer is generated by $M_{\lambda_3} \in \GGp$. We will now 
determine which cones in $\smash{\tSigL}$ are contained in $\sigma_{\lambda_3}$, and are thus by 
Proposition~\ref{prop:stabcont} possibly stabilized by $M_{\lambda_3}$.

\begin{Prop}	\label{prop:twocases}
	If $\tSigL$ contains a cone $\wt{\sigma} \subset \sigma_{\lambda_3}$, which is invariant under the action of 
	$M_{\lambda_3}$, then $\wt{\sigma}$ is of one of the following two forms:
	\begin{enumerate}
		\item	\label{prop:twocases:1st}
			$\wt{\sigma}$ is the $1$--dimensional cone
			\begin{equation}
				\wt{\xi} := \bbR_{\ge 0} \, b_0 \; ,
			\end{equation}
			where $b_0$ is the sum of the generators of the three $1$--dimensional faces of $\sigma_{\lambda_3}$.
		\item	\label{prop:twocases:2nd}
			There is $c_0 \in \sigma_{\lambda_3} \setminus \wt{\xi}$, such that $\wt{\sigma}$ is a 
			$3$--dimensional cone of the form
			\begin{equation}
				\wt{\sigma}(c_0) := \bbR_{\ge 0} \, c_0 + \bbR_{\ge 0} \, M_{\lambda_3} c_0 + 
				\bbR_{\ge 0} \, M_{\lambda_3}^2 c_0 \; .
			\end{equation}
	\end{enumerate}
\end{Prop}
\begin{proof}
	Since $M_{\lambda_3}$ stabilizes the cone $\sigma_{\lambda_3}$, it permutes its three $1$--dimensional faces.
	Furthermore, it has to act transitively on them, since it has order $3$. Hence the only $1$--dimensional cone 
	$\wt{\sigma} \subset \sigma_{\lambda_3}$ which can be left invariant under $M_{\lambda_3}$ is the one given in 
	(\ref{prop:twocases:1st}). Obviously, there can not be a $2$--dimensional cone 
	$\wt{\sigma} \subset \sigma_{\lambda_3}$ which is invariant under the action of $M_{\lambda_3}$, since 
	$M_{\lambda_3}$ has order $3$. For the case of a $3$--dimensional cone, it is easy to see that is has to be of the
	form given in (\ref{prop:twocases:2nd}).
\end{proof}

We can now calculate the singularities of the images of $\orb(\wt{\xi})$ and $\orb(\wt{\sigma}(c_0))$ in $\resptot$
and determine their types.

\begin{Prop}
	\begin{enumerate}
		\item
			If $\wt{\xi} \in \tSigL$, then the image of $\orb(\wt{\xi})$ in $\resptot$ contains exactly three 
			singular points. Two of these are of type $V_{{\frac{1}{3}} (1,1,2)}$, the third is of type 
			$V_{{\frac{1}{3}}(0,1,2)}$. All three are canonical.
		\item
			If $\wt{\sigma}(c_0) \in \tSigL$ for some $c_0 \in \sigma_{\lambda_3} \setminus \wt{\xi}$, then the 
			image of $\orb(\wt{\sigma}(c_0))$ in $\resptot$ consists of one singular point. The singularity is of
			type $V_{{\frac{1}{3}}(0,1,2)}$ and is therefore canonical.
	\end{enumerate}
\end{Prop}
\begin{proof}
	\begin{enumerate}
		\item
  	      	Since we have an explicit description of both $\smash{\wt{\xi}}$ and $M_{\lambda_3}$, we can calculate 
			the singularities directly. The calculation is a straightforward exercise in toric geometry (though
			somewhat lengthy), so we omit it here and give an outline instead. By following the idea of Brasch 
			\cite[Kapitel~3, \S5]{bra:smafpt94}, i.e. passing on to a suitable sublattice, we can describe the 
			action of $M_{\lambda_3}$ on $\smash{\orb(\wt{\xi})}$. It can be described locally near each of its three
			fixed points as a linear operation on $\bbC^3$ (cf. \cite[p.~97]{car:qeaga57}). By calculating the 
			eigenvalues of the matrices representing these linear operations, we can determine the types of the 
			three singularities. Using the Reid--Shepherd--Barron--Tai criterion 
			(cf. \cite[Theorem~4.11]{rei:ypg87}), we can conclude that they are all canonical as claimed.
		\item
			We could proceed similarly to the first case to calculate the action of $M_{\lambda_3}$ on 
			$\orb(\wt{\sigma}(c_0))$. However, there is the following direct argument: Denote the three 
			$1$--dimensional faces of $\wt{\sigma}(c_0)$ by $\tau_i$, $i=1,2,3$. We introduce coordinates 
			$(t_1,t_2,t_3)$ on $T_{\wt{\sigma}(c_0)} \cong \bbC^3$ such that the orbits $\orb(\tau_i)$, $i=1,2,3$,
			are given by
			\begin{equation*}
				\orb(\tau_i) = \left\{ t_i = 0, t_j \neq 0 \text{ for } i \neq j \right\} \; .
			\end{equation*}
			Using the facts that $M_{\lambda_3}$ permutes the three cones $\tau_i$, $i=1,2,3$, and that 
			$M_{\lambda_3}$ has order $3$, we can conclude that the induced action of $M_{\lambda_3}$ on the toric
			variety $T_{\wt{\sigma}(c_0))}$ is of the form
			\begin{equation*}
				\begin{pmatrix}
					t_1 \\ t_2 \\ t_3
				\end{pmatrix}
				\mapsto
				\begin{pmatrix}
					t_3 \\  t_1 \\ t_2 \\
				\end{pmatrix} \; .
			\end{equation*}
			The origin in $\orb(\wt{\sigma}(c_0)) = \left\{ (0,0,0) \right\}$ is fixed by this action. The operation
			is already linear and its eigenvalues are $1$, $\zeta$ and $\zeta^2$, where $\zeta$ denotes a primitive
			third root of unity. Therefore the quotient singularity is of type $V_{{\frac{1}{3}}(0,1,2)}$ as claimed
			and hence canonical.
	\end{enumerate}
\end{proof}

We summarize our discussion in the following

\begin{Lem}	\label{prop:cor2bdc}
	The open corank~$2$ boundary component $\bcD{h}$ of $\resptot$ contains at most canonical singularities.
\end{Lem}

Combining Propositions~\ref{prop:cor1bdc} and \ref{prop:cor2bdc} we obtain the following

\begin{Cor}	\label{cor:canbd}
	The boundary $\resptot \setminus \Ap$ contains at most canonical singularities.
\end{Cor}

With respect to the extendibility of pluricanonical forms, which we are concerned with in the first place, we can state
our results as follows:

\begin{Cor}	\label{cor:obsbd}
	Let $f \otn$ be a pluricanonical form on $\resptot$. Then $f \otn$ extends to a pluricanonical form on any resolution
	of the singularities in $\resptot$.
\end{Cor}
\begin{proof}
	We have seen at the end of the previous section that $\resptot$ contains at most canonical singularities away from
	the boundary. Hence by Corollary~\ref{cor:canbd} $\resptot$ has only canonical singularities. Therefore $f \otn$
 	can be extended as desired.
\end{proof}

\section{Conclusions}	\label{sec:concl}

We can now establish our main result by assembling our results of the previous sections.

\begin{Thm}	\label{thm:gentype}
	$\Ap$ is of general type if there exists a nontrivial weight~$2$ cusp form for $\Gp$ and $p \ge 37$ is prime.
\end{Thm}
\begin{proof}
	Let $f:=f_2^n f_n$ be a weight~$3n$ cusp form for $\Gp$, where $f_2$ is a nontrivial weight~$2$ cusp form and 
	$f_n$ a cusp form of weight $n$. Then $f \otn$ defines a $n$--canonical form on $\Ap$ away from the branch locus.
	The obstructions to extending $f \otn$ to a pluricanonical form on a resolution of singularities of $\smash{\tApk}$ 
	lying in $\Ap$ are given in Lemma~\ref{lem:obsint}. Corollary~\ref{cor:obsbd} shows that these are in fact the only 
	obstructions to extending $f \otn$ to a smooth projective model of $\Ap$. Comparing the dimension of the space of
	our cusp forms $f$ calculated in Proposition~\ref{prop:dimcusp} with these obstructions, we conclude that $\Ap$ is
	of general type as long as the condition
	\begin{equation*}
		\frac{p^2 + 1}{8640} > \frac{7}{54} - \frac{1}{3p}
	\end{equation*}
	is satisfied. This is the case, if $p \ge 37$.
\end{proof}

Using Proposition~\ref{prop:nontrivcusp} we obtain the following immediate

\begin{Cor}
	\begin{enumerate}
		\item
			$\Ap$ is of general type for $p > 71$.
		\item
			$\Ap$ is of general type for $p = 37,43,53,61$ and $67$.
	\end{enumerate}
\end{Cor}

It is rather unlikely that the given bound is sharp. In general, there will be weight~$3n$ cusp forms $f$ which can not
be expressed in the form $f_2^n f_n$ but which extend nevertheless to pluricanonical forms on a smooth projective model
of $\Ap$. The only place where we made use of the special form of our cusp forms was in Lemma~\ref{lem:ext} to guarantee
the extendibility to the smooth part of the boundary. Hence, if one starts with a different space of weight~$3n$ cusp 
forms which extend to the smooth points in the boundary, one should be able to proceed analogously to extend them to the
whole of a smooth projective model of $\Ap$.

There should be an analogous result for the space $\Apol{t}$, where $t$ does not need to be prime any more. In this case,
however, the singularities tend to be much more complicated and so the obstructions will be much larger. If $t$ is of a 
special form, though, we can conclude that $\Apol{t}$ is of general type using Theorem~\ref{thm:gentype}.

\begin{Cor}
	If $q \in \bbN$ and $p > 71$ or $p \in \left\{ 37,43,53,61,67 \right\}$ is prime then the space $\Apol{pq^2}$ is
	of general type.
\end{Cor}
\begin{proof}
	According to \cite[Corollary~7.2]{san:mpas97} there is a surjective morphism $\Apol{pq^2} \rightarrow \Ap$ which
	can be extended to suitable smooth compactifications. Hence $\Apol{pq^2}$ is of general type as long as $\Ap$ is.
\end{proof}

\vspace{1cm}

\begin{tabular}{l}
Cord Erdenberger\\
Universit\"at Hannover\\
Institut f\"ur Mathematik\\
Welfengarten 1\\
D-30167 Hannover\\
Germany\\
\\
{\tt erdenber@math.uni-hannover.de}
\end{tabular}

\end{document}